\documentclass{article}

\usepackage{arxiv}

\usepackage[utf8]{inputenc} 
\usepackage[T1]{fontenc}    
\usepackage{hyperref}       
\usepackage{url}            
\usepackage{booktabs}       
\usepackage{amsfonts}       
\usepackage{nicefrac}       
\usepackage{microtype}      
\usepackage{lipsum}		
\usepackage{graphicx}

\usepackage{moreverb}

\usepackage{graphicx, subfigure}

\newcommand\BibTeX{{\rmfamily B\kern-.05em \textsc{i\kern-.025em b}\kern-.08em
T\kern-.1667em\lower.7ex\hbox{E}\kern-.125emX}}

\usepackage{graphicx}
\RequirePackage{fix-cm}

\def\beq{\begin{equation}}
\def\eeq{\end{equation}}
\def\baq{\begin{eqnarray}}
\def\eaq{\end{eqnarray}}
\def\bc{\begin{center}}
\def\ec{\end{center}}
\def\ds{\displaystyle}

\def\eps{\varepsilon}

\def\gr{\gamma_{r}}
\def\gr1{\gamma_{r1}}

\usepackage{hyperref}
 \usepackage{epsfig} 
 \usepackage{amssymb}
 \usepackage{amsmath}
\usepackage{multirow}
\usepackage{epsfig} 
\input{epsf}
\usepackage{graphicx}
\usepackage{graphicx, subfigure}
\usepackage[utf8]{inputenc}
\usepackage[english]{babel}

\def\beq{\begin{equation}}
\def\eeq{\end{equation}}
\def\baq{\begin{eqnarray}}
\def\eaq{\end{eqnarray}}
\def\bal{\begin{align} }
\def\eal{\end{align} }
\def\bc{\begin{center}}
\def\ec{\end{center}}
\def\ds{\displaystyle}

\def\eps{\varepsilon}

\def\gr{\gamma_{r}}
\def\gr1{\gamma_{r1}}

\def\tt{\theta}

\def\eps{\varepsilon}

\def\gr{\gamma_{r}}
\def\gr1{\gamma_{r1}}

\def\eps{\varepsilon}

\def\til{\tilde{\theta}}

\title{Convergence Rate Improvement of Richardson and Newton-Schulz Iterations}

\author{ Alexander Stotsky \\
   Department of Computer Science and Engineering \\
     Chalmers University of Technology  \\
     Gothenburg  SE - 412 96, Sweden  \\
	\texttt{alexander.stotsky@chalmers.se} \\
        \texttt{alexander.stotsky@telia.com}  }

\hypersetup{
pdftitle={A template for the arxiv style},
pdfsubject={q-bio.NC, q-bio.QM},
pdfauthor={David S.~Hippocampus, Elias D.~Striatum},
pdfkeywords={First keyword, Second keyword, More},
}

\begin{document}
\maketitle

\begin{abstract}
~Fast convergent, accurate, computationally efficient, parallelizable,
and robust matrix inversion and parameter estimation algorithms are required
in many time-critical and accuracy-critical applications
such as system identification, signal and image processing, network and big data analysis, machine learning and in many others.
\\
This paper introduces  new composite power series expansion
with optionally chosen rates
(which can be calculated simultaneously on parallel units with different computational capacities)
for further convergence rate improvement
of high order Newton-Schulz iteration.
New expansion was integrated into the Richardson iteration and resulted in significant
convergence rate  improvement.
The improvement is quantified via explicit transient models
for estimation errors and by simulations.
In addition, the recursive and computationally efficient version of the combination of
Richardson iteration and Newton-Schulz iteration with composite expansion is developed for
simultaneous calculations.
\\ Moreover, unified factorization is developed in this paper in the form
of tool-kit for power series expansion, which results in  a new family
of computationally efficient Newton-Schulz algorithms.
\end{abstract}

\maketitle

\section{Introduction}
\noindent
Least squares method is widely used in control, system identification, signal processing,
\cite{lju1} - \cite{gus1}, statistics, \cite{bat1} as well as in  many computational applications such as emerging big data applications, \cite{wol}, machine learning, \cite{lag} and in many other areas.
For accurate solution many least squares problems
(for example the problems related to data, signal and image processing,
system identification, network analysis and many others) can be associated
with calculation of the parameter vector $\tt_*$, which
satisfies the algebraic equation
\beq
A \tt_* = b \label{ae1}
\eeq
where $b$ is the vector, and the matrix $A$ is SPD (Symmetric and Positive Definite)
matrix. For example, the matrix $A$ is SPD for the systems with harmonic regressor,
\cite{fom}, \cite{bay1} and multiplication of any invertible matrix $A$ by its transpose transforms the system to the SPD case with the Gram matrix, \cite{bjor}.
The numerical stability problems associated with ill-conditioning of the Gram matrix
can be solved using different types of preconditioning techniques, see for example,
\cite{ben}, \cite{st15} and references therein (see also Section~\ref{comp} for
simulations of the ill-conditioned matrices).
\\ Iterative methods for solving (\ref{ae1}) are often preferable (especially for large-scale
systems) due to simplicity, better accuracy and robustness,
less processor time and memory space compared to direct methods.
The most general and well-known method for iterative calculation of the matrix inverse
is high order Newton-Schulz algorithm described in \cite{isa} - \cite{pan}
and in many other publications.
The second order version of Newton-Schulz iteration, see
for example \cite{sch} - \cite{beniz} is the most known.
\\ High order Newton-Schulz algorithms are well-discussed
in the literature. However, the questions associated with
the relation between high order Newton-Schulz algorithms and power series expansions
were not properly  studied. The paper \cite{jan},
which was the first paper with the description of the relation between second order Newton-Schulz algorithm and power series expansion does not provide  complete description of this relation.
\\ Reduction of the computational complexity of high order Newton-Schulz algorithm
is one the most important challenges in this area. The computational complexity can be reduced
via factorizations of the power series, see for example \cite{sti}, \cite{pan},
\cite{sha} - \cite{bur} and references therein.
Practical applications of these factorizations (excepting Horner's rule)
are hampered by the lack of general unified description.
\\ Computational resources with high degree of parallelism
(instead of single computing units) will be available
in the future for implementation of numerical methods.
The computational performance of iterative solvers can also be increased
via parallel computing (especially for large scale systems),
achieved for example, via multiprocessor and virtual systems, \cite{saa1} - \cite{pen}.
In order to improve the performance a serial algorithm is usually
converted to  parallel algorithm, see for example \cite{van}.
 This paper proposes a new approach for convergence rate improvement
where novel iterative algorithms are designed with high degree of parallelism
(or enhanced parallelism).
In other words, the iterative algorithm is designed as a number of independent
computational parts (the number of parts is associated with the degree of
parallelism) which can be executed simultaneously.
The challenges associated with computational efficiency are addressed already on the design level in this case,
providing new opportunities for high performance parallel processing.
\\ This paper introduces  new composite power series expansion
with optionally chosen rates and high degree of parallelism  for further convergence rate improvement in the unified framework described in \cite{st19}. New expansion applied to
Richardson iteration resulted in significant improvement of the convergence rate.
Simulation results are presented for quantification of the
improvements of new algorithms compared to recent algorithms
described in \cite{st19}. Moreover, explicit transient
models are derived for all new algorithms described in this paper.
In addition, the recursive and computationally efficient version of the combination of
Richardson iteration and Newton-Schulz iteration with composite expansion is developed for
simultaneous calculations.
\\ Finally, factorization tool-kit is developed in
this paper for general power series expansion, which
allows nested applications and results in a family
of new computationally efficient algorithms.
\\ This paper is organized as follows. The paper starts with the representation of
Newton-Schulz  iteration as power series expansion
in Section~\ref{nssec1}. A unified power series factorization for reduction of computational complexity
is introduced in  Section~\ref{psf} and represented in the form of tool-kit
in Section~\ref{appen}.
 New high order Newton-Schulz algorithms with
composite polynomial are presented in Section~\ref{cns}.
Richardson iteration with high order convergence accelerator is described
in Section~\ref{ra} and compared to existing algorithms by simulation in Section~\ref{comp}.
The paper ends with brief conclusions in Section~\ref{conc}.
\\ \\ {\it This paper was presented on the  21-st IFAC World Congress in Berlin, Germany, July 12-17, 2020,
\cite{ifac2020}.}

\section{Splitting \& Preconditioning}
\label{splpre}
\noindent
Numerical solution of the system of linear equations (\ref{ae1})
using power series expansions requires splitting and preconditioning.
Any positive definite and symmetric matrix $A$, whose inverse should be calculated can be split as follows,
see for example \cite{che} and references therein :
\baq
 & & A = S - D \label{minus} \\
 & & I - S^{-1} A =  S^{-1} D \label{sp1} \\
 & & \rho(I - S^{-1} A) =  \rho(S^{-1} D) < 1 \label{spro}
\eaq
where the spectral radius $\rho(\cdot)$ defined in (\ref{spro}) is less than one
for symmetric and positive definite matrices $A$ and $S$ (where $S^{-1}$ is the preconditioner),
provided that $2 S - A$ is a  positive definite matrix, \cite{hac}.
\\ For example, the matrix $S$  can be chosen as a diagonal matrix,
which contains the diagonal elements of SDD (Strictly Diagonally Dominant) and
positive definite  matrix $A$, see \cite{hor} for the general case
and \cite{st10} for systems with harmonic regressor.
\\ For positive definite (not SDD) matrix $A$ the simplest preconditioner
can be chosen as  $\ds S^{-1} = I / \alpha$ with
$\alpha = \| A \|_\infty / 2 + \eps$, where $ \| \cdot \|_\infty $ is the maximum row sum matrix norm,
and $\eps > 0$ is a small positive number, \cite{st15}, \cite{sha1}.
\\ Other types of preconditioning can be found in \cite{che}, \cite{hac}, \cite{var}, \cite{st161},
see also references therein.

\section{Newton-Schulz Iteration as Fast Power Series Expansion }
\label{nssec1}
\noindent
 The results presented in this Section introduce
computationally efficient factorization of initial power series
and show several steps (step by step) of fast matrix power
series expansion which coincide with  Newton-Schulz iteration.
The relation between Newton-Schulz approach and
power series expansion opens new opportunities
for reduction of the computational complexity of
Newton-Schulz algorithms.
\\ The following initial power series factorization :
\baq
 G_0 &=& \sum_{j = 0}^{w-1} (S^{-1}D)^{(p+1) j} ~\{ \sum_{d=0}^{p} (S^{-1}D)^{d} \} S^{-1}
 \label{g01} \\
  &=& \sum_{j=0}^{h-1} (S^{-1}D)^{j} S^{-1} = (I -   (S^{-1}D)^{h}) A^{-1} \label{g02} \\
  F_0 &=& I -  G_0 A = (S^{-1}D)^{h} \label{g03}
\eaq
where (\ref{g03}) defines initial inversion error, and  $p = 0,1,2,...$,
$w = 1,2,3,...$, and $h = w (p+1) = 1,2,3 ... $,
gives the starting point  for the following steps of Newton-Schulz iteration:
\par \underline{Step 1.}
\baq
& & G_1 = \{ \sum_{j=0}^{n-1}  F^j_0 \} ~ G_0 = \{ \sum_{j=0}^{n-1} (S^{-1}D)^{h j} \} ~ G_0  \nonumber \\
& & [ I + (S^{-1}D)^{h} + (S^{-1}D)^{2 h} + ... +  (S^{-1}D)^{(n-1) h} ] \nonumber \\
& & [ I + (S^{-1}D) + (S^{-1}D)^{2} + ... +  (S^{-1}D)^{h-1} ]  S^{-1} \nonumber \\
 & & = [ I + (S^{-1}D)  + ... +  (S^{-1}D)^{(h n - 1) } ]  S^{-1} \label{ser1}  \\
 & & F_1 = I -  G_1 A = (S^{-1}D)^{h n} \label{f10}
\eaq
where $G_1$ in (\ref{f10}) is calculated via (\ref{ser1}).
\par \underline{Step 2.}
\baq
& & G_2 = \{ \sum_{j=0}^{n - 1}  F^j_1 \} ~ G_1 = \{ \sum_{j=0}^{n-1} (S^{-1}D)^{h n j} \} ~ G_1  \nonumber \\
 & & F_2 = I -  G_2 A = (S^{-1}D)^{h n^2} \label{f20}
\eaq
Further evaluation in  \underline{Step k} gives classical high order Newton-Schulz algorithm (\ref{ns}) and error model
(\ref{fk1}) :
\baq
& & G_k = \{ \sum_{j=0}^{n-1}  F^j_{k-1} \} ~ G_{k-1} \label{ns} \\
 & & F_k = I -  G_k A = F_0^{n^k}  = (S^{-1}D)^{h n^k} \label{fk1}
\eaq
where  $G_k$ is estimate of $A^{-1}$, $n = 2,3, ...$ and  $k = 1,2,3,...$
\\ Notice that the factorization similar to  (\ref{g01}) can be applied to the power series
(\ref{ns}) for improvement of computational efficiency.
To this end the unified factorization method is developed in the next Section.

\section{Reduction of Computational Complexity via Unified Factorization: Nested Algorithms}
\label{psf}
\begin{figure}
\centerline{\psfig{figure=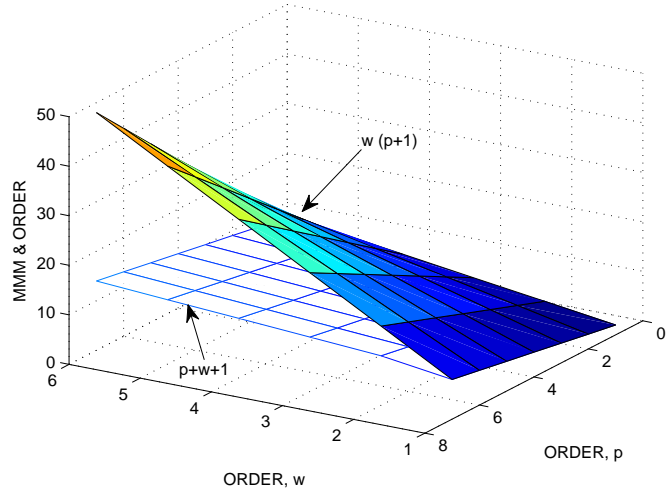,height=70mm}}
\begin{center}
\caption{
The order $ h = w ~ ( p + 1) $ (which is equal to the number of mmm  for conventional realization), where $p = 1,2,...,7$, $~ w = 1,2,...,6 $ of the factorization (\ref{xa2})-(\ref{z2}) is plotted with colored surface and the number of mmm  $ N_p = p + w + 1$  is plotted with a white surface. }
\label{mmm1}
\end{center}
\end{figure}

\noindent
Consider the following matrix power series:
\baq
   Z &=& \{ \sum_{j=0}^{h-1} Y^j \} ~ X
   \label{z1} \\
   Y &=& I - X A \label{xa1}
\eaq
where $X,Y,Z$ are matrices of corresponding dimensions, $I$ is the identity matrix,
$h = 2,3,4, ... $~. Realization of the algorithm (\ref{z1}) requires $h$ mmm (matrix-by-matrix
multiplications) per iteration loop according to Horner's scheme,
see for example \cite{pan}, \cite{st19}.
\\
Notice that Horner's rule is not optimal for evaluating matrix polynomials, \cite{sti}
and for reduction of the computational complexity the power series (\ref{z1})
can be factorized as follows:
\baq
    Y &=& I - X A \label{xa2} \\
  U &=& \{ \sum_{d=0}^{p} Y^d \} ~ X
  \label{u1} \\
  Y^{p+1} &=& I -  U A \label{ua1} \\
  Z &=& \{ \sum_{j=0}^{w-1} Y^{(p+1) j} \} ~ U    \label{z2} \\
   h &=& w ~ ( p + 1), ~~ p = 0,1,2, ..., ~ w = 1,2,3, ...  \label{h1} \\
   N_p &=& p + w + 1   \label{np1}
\eaq
where $N_p$ is the number of multiplications
for realization of the algorithm (\ref{xa2}) - (\ref{z2}) of the order $h$ defined
in (\ref{h1}). Indeed, realization of (\ref{xa2}) - (\ref{ua1}) requires $p + 2$ multiplications, and
realization of the series (\ref{z2}) which can be calculated as follows:
\beq
Z_i = Y^{p+1} Z_{i-1} + U, ~ for~ i = 1 : (w-1),~ Z_0 = U
\label{zi1}
\eeq
requires $w-1$ multiplications.
\\ Notice that $N_p = p + 1$ and  $N_p = w $ in (\ref{np1}) for the case
where  $w = 1$ and $p=0$ respectively.
\\ Notice that the idea of factorization (\ref{xa2}) - (\ref{z2})
is associated with the Newton-Schulz iteration (\ref{ns}), where
the sum $\ds \sum_{j=0}^{w-1} Y^{(p+1) j}$ in (\ref{z2})  corresponds to
$\ds \sum_{j=0}^{n-1}  F^j_{k-1}$ and $U$ in (\ref{u1}) is associated with $G_{k-1}$.
\\ Representation (\ref{z1}) and (\ref{xa1}) includes
Newton-Schulz algorithm (\ref{ns}), (\ref{fk1}) with $X = G_{k-1}$, $Y = F_{k-1}$,
and $Z = G_{k}$. In addition, equation (\ref{z2}) represents (\ref{g01}) with
$ Y = (S^{-1}D)$ and $X = S^{-1}$ and can be used for computationally efficient calculations
of the initial power series. Algorithm (\ref{xa2})-(\ref{z2}) describes unified and
systematic way for power series factorization, order reduction and improvement
of the computational efficiency. Application of the algorithm (\ref{xa2}) - (\ref{z2})
to factorization of Newton-Schulz iteration of orders $2 - 19$ is
demonstrated in the form of tool-kit in Appendix, see Section~\ref{appen}.
\\ The number of mmm for conventional recursive realization of high order Newton-Schulz algorithm
is equal to the algorithm order. Factorization (\ref{xa2})-(\ref{z2}) reduces the number of mmm
to (\ref{np1}) for the order (\ref{h1}). The reduction of computational complexity is quantified
in Figure~\ref{mmm1}, where the order (which is equal to the number of mmm for
conventional realization) is plotted with colored surface
and the number of mmm for (\ref{xa2})-(\ref{z2}) is plotted with a white surface.
The complexity can be essentially reduced for higher orders.
The efficiency index introduced in \cite{ehr} has the following form for the factorization (\ref{xa2})-(\ref{z2}):
\beq
 EI = \ds  [w ~ ( p + 1)]^{1 / (p + w + 1) } \label{ei}
\eeq
Notice that sequential application of the factorization (\ref{xa2})-(\ref{z2}) implies
further reduction of the computational complexity for higher orders.
The sum  (\ref{z2}) can also be calculated more efficiently for specific orders
compared to (\ref{zi1}). Nested application of (\ref{xa2})-(\ref{z2}) is
illustrated by the following example.
\\
\underline{Example:} Nested algorithm as unification of the hyperpower iteration method
(described in \cite{bur}) of order $45$ that requires $10$ mmm only.
\\
Newton-Schulz
iteration of the order $45$ can be factorized in a step-wise way
as follows:
\begin{align}
& Z = (I + Y^9 + ... + Y^{36}) ~~ [ I + Y + ... + Y^8] ~ X  \label{ex1} \\
 &= (I + Y^9 + ... + (Y^9)^4) ~~ [I + Y^3 + Y^6]~[I + Y + Y^2] ~ X  \label{ex2} \\
 & = \{ I + ( I + (Y^9)^2 ) ~ ( Y^9 +(Y^9)^2) \} ~ ~ [I + Y^3 + Y^6]  \nonumber \\
 & [3I + (XA)~(- 3I + XA)] ~ X   \label{ex3} \\
 &  Y^3 =  I -  [\sum_{d=0}^{2} Y^d] ~X~ A, ~  Y^9 =  I - [\sum_{d=0}^{8} Y^d] ~X~ A  \nonumber
 \end{align}
Algorithm (\ref{ex1}) represents the algorithm (\ref{xa2}) - (\ref{np1}) of the order $h=45$ defined
in (\ref{h1}) with $p=8$ and $w = 5$ and requires $ p + w + 1 = 14 $ mmm.
Further application of the algorithm (\ref{xa2}) - (\ref{np1}) to the eighth order polynomial\footnote{
Other type of factorization of the eighth order polynomial is presented  \\
in Appendix, Table~\ref{t2} for h = 9}
 in (\ref{ex1}) and factorization of the fourth order polynomial (w.r.t. $Y^9$)
results in algorithms (\ref{ex2}), (\ref{ex3}), which requires ten mmm only with EI $ =  1.4633 $.
\\
Notice that nested method for derivation of the
computationally efficient algorithms
based on (\ref{xa2}) - (\ref{np1}) is more simple compared to the method
described in \cite{bur}. Moreover,  the method is universal
(compared to heuristic methods) and applicable to any order,
see Section~\ref{appen}.
\\
The efficiency index EI does not account for robustness
with respect to error accumulation. Minimization of the number of mmm makes
iteration more robust.
\\ Notice also that the efficiency index  of classical Newton-Schulz
iteration of the second order, which is the most robust,  is EI $ = 1.4142 $,
see for example \cite{sol}, \cite{solor} for comparisons of the efficiency indexes.

\section{Novel Newton-Schulz Algorithm with Composite Polynomial
and Enhanced Parallelism for Simultaneous Calculations}
\label{cns}
\noindent
The unified framework for convergence rate improvement of high order Newton-Schulz matrix inversion algorithms was proposed in \cite{st19} .
The following new composite power series expansion for
Newton-Schulz iteration with different expansion rates for further convergence rate improvement extends this framework as follows :

\baq
G_k &=& \underbrace{T_c}_{\begin{subarray}{l}\text{Composite }\\
    \text{Polynomial }\end{subarray}}
 + \underbrace{\Gamma_c}_{\begin{subarray}{l}\text{Composite }\\
    \text{Residual}\end{subarray}} ~~ \underbrace{\{ \sum_{j=0}^{n-1}  F^j_{k-1} \} ~ G_{k-1}}_
 {\begin{subarray}{l}\text{Newton-Schulz }\\
    \text{Iteration }\end{subarray}} \label{gkgen1} \\
T_c &=& \sum_{i=1}^{w} \{ \prod_{p = 0}^{i-1} \Gamma_p \}  T_i = T_1 + \Gamma_1 T_2 + ... + \prod_{p=1}^{w-1} \Gamma_p T_w  \label{tc} \\
 \Gamma_c &=& \prod_{p=1}^{w} \Gamma_p \label{gc}
\eaq
where $\ds T_c$ is the composite power series expansion and the
composite residual is defined as the product of
the residual terms $\ds \Gamma_c = \Gamma_1 \Gamma_2 \Gamma_3 ... \Gamma_w$
with the spectral radius $\rho(\Gamma_i) < 1$, and $\Gamma_0 = I$.
Power series expansions $T_i$ satisfy the following relations:
\baq
   T_i A &=& I - \Gamma_i~,~~~ i = 1,...,w \label{rterm1} \\
   T_c A  &=& I - \Gamma_c \label{rel1}
\eaq
Multiplication of both sides of equation (\ref{gkgen1}) by $A$ together with the relation (\ref{rel1})
results in the following error model:
\beq
  F_k = \Gamma_c F^n_{k-1} \label{gem1}
\eeq
where the following spectral radius is less than one, $\rho(\Gamma_c) \le \rho(\Gamma_1) ...  \rho(\Gamma_w) < 1$
according to Gelfand's formula provided that the matrices $\Gamma_i$ commute.
\\ The error model (\ref{gem1}) with composite power series expansion is the same as the error model (22) in
\cite{st19}  for a single power series. The advantages of composite expansion are discussed below.
\\ The power series expansions $T_i$ (which can be calculated simultaneously
on parallel computational units) can be taken as
\baq
  T_i &=& \{ \sum_{j=0}^{x_i - 1} (S^{-1} D)^j  \} S^{-1} = (I - \Gamma_i) A^{-1} \label{ti} \\
  \Gamma_i &=& (S^{-1} D)^{x_i} \label{ga1}
\eaq
The rate of expansion $x_i$ can be chosen using computational capacity of each parallel computational
unit (fast power series expansion should be implemented on more powerful computational unit).
For example $x_i$ can be taken as a polynomial which is  a function of step number $k$ or as
rapidly expanding power series associated with high order Newton-Schulz iteration with
 $x_i =  m^k$,~$ m = 2,3,... $, see for example \cite{st19} for this and other choices.
\\ Notice that the algorithm (\ref{gkgen1}) - (\ref{gc})  has especially simple
 form for $n = w = m = 1$, $\ds G_k = (S^{-1} D) G_{k-1} + S^{-1}$ which was derived in \cite{bre}
 directly from splitting.

\subsection{Double Newton-Schulz Algorithm with High Order Residual as Convergence Accelerator}
\label{dns}
\noindent
The advantages of the framework described above are especially pronounced for case when choosing
a number of the same rapid expansions $T_i$ (high order Newton-Schulz iterations for example)
with the expansion rate associated with the order $n$ in (\ref{gkgen1}). Fast and computationally efficient
algorithms can be designed in this case.
\\ Consider algorithm (\ref{gkgen1}) with $T_1 = T_2 = ... = T_k$, where $T_k$ and $\Gamma_k$ are defined
in (\ref{ti}),(\ref{ga1}) with $w=n$ and $x_i =  h n^k $, $h,n = 1,2, ...$ :
\begin{align}
\Gamma_k &= I - L_{k-1} A \label{gamk11} \\
L_k &=  \{ \sum_{j=0}^{n-1} \Gamma^j_k \} ~ L_{k-1} \label{lkk1} \\
\Gamma^n_k &=  I - L_k A \label{gamk21}  \\
G_k &= \underbrace{L_k}_{\begin{subarray}{l}\text{Newton-Schulz }\\
    \text{Iteration }\end{subarray}} + \underbrace{\Gamma^n_k}_{\begin{subarray}{l}\text{High Order }\\
    \text{Convergence} \\
    \text{Accelerator} \end{subarray}}   ~  ~ \underbrace{ \{ \sum_{j=0}^{n-1}  F^j_{k-1} \} ~ G_{k-1}}_
 {\begin{subarray}{l}\text{Newton-Schulz }\\
    \text{Iteration }\end{subarray}}
 \label{gkgen21} \\
 F_k &= I -  G_k A = \Gamma^n_k ~ F^n_{k-1}      \label{fks1} \\
F_k &= (S^{-1} D)^{\ds h ~ ( k ~ n^{k+1} + ~ n^k )} \label{ermoex1}
\end{align}
where $\ds L_0 =\{ \sum_{j=0}^{n-1} \Gamma^j_0 \} ~ T_0 $,
$\Gamma_0 = I - T_0 A$, and  $ T_0 = G_0 $ ($G_0$ is calculated via (\ref{g01})) are precalculated.
The algorithm (\ref{gamk11}) - (\ref{ermoex1}) has two Newton-Schulz loops
(which can be calculated simultaneously) of the same order $n$
associated with inversion errors (\ref{gamk11}) and (\ref{fks1}). The sums
 $ \ds \sum_{j=0}^{n-1} \Gamma^j_{k} $   and $ \ds \sum_{j=0}^{n-1}  F^j_{k-1} $
 can be calculated recursively using Horner's scheme, \cite{st19}
 or factorizations, see Section~\ref{appen}.
 \\ Notice that the algorithm derived in \cite{bre} and the algorithm (15) in
 \cite{st19} are special cases of the  algorithm (\ref{gamk11}) - (\ref{gkgen21})
 for $n=h=1$ and $n=1$ and $h \ge 1$ respectively.
  \\ \underline{Remark 1.} Comparison of the error model (\ref{ermoex1}) with the error model
(\ref{fk1}) of classical high order Newton-Schulz algorithm shows that
the algorithm (\ref{gamk11}) - (\ref{fks1}) has significantly higher convergence rate
due to the term  $\ds k ~ n^{k+1}$.
\\ \underline{Remark 2.} The algorithm similar to (\ref{gamk11}) - (\ref{ermoex1}) was proposed in \cite{sri}.   The algorithm written in the following form:
\baq
 Z _{k} &=&  \sum_{j = 0}^{p-1} (I - Z_{k-1} A)^j  ~ Z_{k-1} \label{sri1} \\
  G _{k} &=&  G_{k-1} + (I - G_{k-1} A)  ~ Z_{k} \label{sri2}
\eaq
has also two Newton-Schulz loops, where both $Z_k$ and $G_k$ are the estimates
of the matrix inverse and $p= 2,3,...$ is the order.
\\
Algorithm (\ref{sri1}), (\ref{sri2}) has the following error model
\baq
L_k &=& L_{k-1}^p, ~~ L_k = I - Z_k A       \label{srilk}   \\
F_k &=& F_{k-1} L_{k-1}^p, ~~ F_k = I - G_k A   \label{srifk1}
\eaq
where $L_k$ and $F_k$ are estimation errors.
\\ The algorithm (\ref{gamk11}) - (\ref{ermoex1}) has faster convergence
due to the high order error $F_{k-1}^n$  in the error model
(\ref{fks1}) compared to algorithm  (\ref{sri1}), (\ref{sri2}) which has the error model (\ref{srifk1}) with the first order error $F_{k-1}$.

\section{Richardson Iteration with High Order Convergence Accelerator}
\label{ra}
\noindent
\subsection{Algorithm Description}
\noindent
Combination of Richardson iteration, see \cite{ric} and matrix inversion algorithms was proposed first
in \cite{dub} for improvement of the convergence rate of estimated parameters.
A number of combinations of Richardson iteration with matrix inversion techniques
has been developed in recent years, see for example \cite{st15} and \cite{che1} - \cite{sri1}
and references therein. Unified framework for many combinations was
proposed recently in \cite{st19}.
New matrix inversion algorithms described in the previous Section can be integrated
into the  Richardson iteration within this unified framework.
\\ The parameter vector
in (\ref{ae1}) can be estimated via recursive algorithm as follows:
\beq
\tt_k = \tt_{k-1} -
\underbrace{
[ ~ L_k + \Gamma_k^n ~ \{ \sum_{j=0}^{q-1} F^j_k \} ~ G_k ]}_{\begin{subarray}{l}\text{Fast Matrix Inversion }\\
    \text{Algorithm}\end{subarray}}
 ~ \underbrace{\{A \tt_{k-1} - b \}}_{\begin{subarray}{l}\text{Parameter}\\
    \text{Estimation Error}\end{subarray}}
\label{ttk1}
\eeq
where $\tt_{k}$ is the estimate of $\tt_*$  and
$\Gamma_k$, $L_k$, $F_k$ and $G_k$ are calculated in (\ref{gamk11}) - (\ref{fks1})
and $q=1,2,...$ is the order of Neumann series.
\\ The following model is valid for estimation error $\til_k = \tt_k - \tt_*$ :
\baq
 \til_k &=& \Gamma_k^n ~ F^q_k  ~ \til_{k-1}  \label{riem1}  \\
  \til_k &=& (S^{-1} D)^{\ds h ~ \{ ( k ~ n^{k+1} + ~ n^k ) ~ q + n^{k+1} \} } ~ \til_{k-1}                                              \label{riem2}
\eaq
The error model has especially simple form for the case where $q=n$:
\begin{align}
 \til_k &= \Gamma_k^n ~ F^n_k  ~ \til_{k-1}  \label{riem3}  \\
  \til_k &= (S^{-1} D)^{\ds h ~ ( k ~ n^{k+2} + ~ 2~ n^{k+1})  } ~ \til_{k-1}                                              \label{riem4} \\
  \til_k &= (S^{-1} D)^{\ds \gamma_k} ~ \til_0  \label{riem5} \\
  \gamma_k &= h ~ n^2~ \frac{(k ~ n^{k+2} -(k-1)~ n^{k+1} - 2 ~ n^{k} - n + 2)}{(n-1)^2}
  \label{riem6}
\end{align}
where $n > 1$, $\til_0 = \tt_0 - \tt_*$ and  $ \tt_0 = L_0 b $.
\\ The error model (\ref{riem5}) shows significant improvement of the convergence rate
of estimated parameters in the algorithm (\ref{ttk1}). This improvement is associated with
introduction of the fast matrix inversion algorithms in the Richardson loop,
and it is quantified in the next Section.
\\ Notice that the algorithm (\ref{ttk1}) can also be seen as an extension
of the unified framework of Richardson iteration, \cite{st19}
where the multiplicative high order accelerator $\Gamma_k^n$
and additional Newton-Schulz loop  $L_k$ were introduced for
convergence rate improvement.
\\ \underline{Remark 3.} Stability analysis of combinations of Richardson iteration
and matrix inversion methods described in \cite{dub} and
\cite{che1},  \cite{sri1} is based on the residual error model
$ r_k = b - A \tt_k$, whereas the analysis in
\cite{st15}, \cite{st19} and \cite{st14} (including the analysis above)
is presented in terms of the parameter mismatch $\til_k$, where
$ A^{-1} r_k = \tt_* - \tt_k = - \til_k$.
Notice that the parameter mismatch is widely used in the area of
system identification for  stability analysis, \cite{lju1} - \cite{gus1}
and allows simplified representation of the error models in
unified Richardson and Newton-Schulz framework.
Such error models simplify essentially the stability analysis, which allows
integration of more sophisticated algorithms (which in turn could essentially
improve convergence rate) into the framework.
\subsection{Reduction of Computational Complexity via Recursive and Simultaneous
Calculations}
\noindent
For development of the computationally efficient version the algorithm (\ref{ttk1})
is presented in the following form:
\baq
\tt_k &=& \tt_{k-1} - \omega_k~(A \tt_{k-1} - b) \label{redttk1} \\
\omega_{k} &=& L_k + \Gamma^n_k  ~ \{ \sum_{j=0}^{n - 1} F^j_k \} ~ G_k
\label{redttk2}
\eaq
Recursive algorithm for calculation of  $\omega_{k}$ described below
is divided in two independent computational parts for simultaneous calculations.
\\
Calculations of both parts start with calculation of the
$\ds  \sum_{j=0}^{n-1} \Gamma^j_k $, where
$\ds \Gamma_k = \Gamma^n_{k-1}$.
\\ 1) The first part is associated with the calculation
of $G_k$ in  (\ref{gkgen21})  via $\omega_{k-1}$ as follows:
\beq
 G_k = [ \sum_{j=0}^{n-1} \Gamma^j_k] ~ [ \omega_{k-1} -
 \sum_{j=0}^{n-1}  F^j_{k-1}  ~ G_{k-1} ] +  \sum_{j=0}^{n-1}  F^j_{k-1}  ~ G_{k-1}
\label{gkrecu1}
\eeq
which requires one matrix multiplication only and further calculation of
$\ds \sum_{j=0}^{n-1}  F^j_{k} ~ G_{k} $ with  $G_k$ defined in (\ref{gkrecu1})
which in turn can be  further divided
in independent parts (and calculated for example according to Horner's scheme or factorizations,
see Section~\ref{appen} ).
\\ 2) The second part is associated with calculations of $L_k$
and $\Gamma^n_k$  in (\ref{lkk1}) and (\ref{gamk21}) respectively
using $\ds  \sum_{j=0}^{n-1} \Gamma^j_k $.
\\ The results of both parts are merged in (\ref{redttk2})
to be included in the Richardson iteration  (\ref{redttk1}).
\\ Notice that the matrix-by-vector product $ A \tt_{k-1} $ in (\ref{redttk1}) can be
   easily calculated in parallel via methods described for example in \cite{saa1}.

\section{Comparisons \& Quantification of the Performance }
\label{comp}
\begin{figure}
\centerline{\psfig{figure=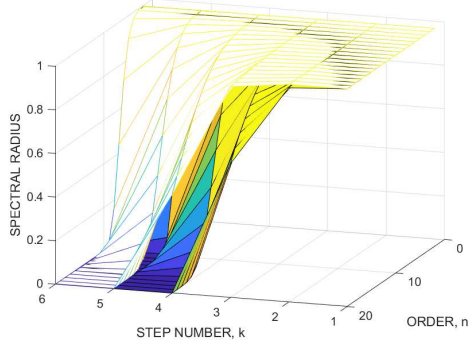,height=50mm}}
\begin{center}
\caption{\small{The Figure shows
the spectral radius  $ \ds \rho^{\ds ( k ~ n^{k+1} + ~ n^k )} $
for double Newton-Schulz matrix inversion algorithm defined in (\ref{ermoex1})
(described in Section~\ref{dns}), plotted as colored surface.
The spectral radius for Newton-Schulz iteration with improved convergence rate
described in \cite{st19},  $ \ds \rho^{\ds  (k +1) ~ n^{k} } $
is  plotted as white surface.
Both surfaces are plotted for  the spectral radius of ill-conditioned case
as functions of the order $n$ and step number $k$.
 }}
\label{figsurface3}
\end{center}
\end{figure}

\begin{figure}
\centerline{\psfig{figure=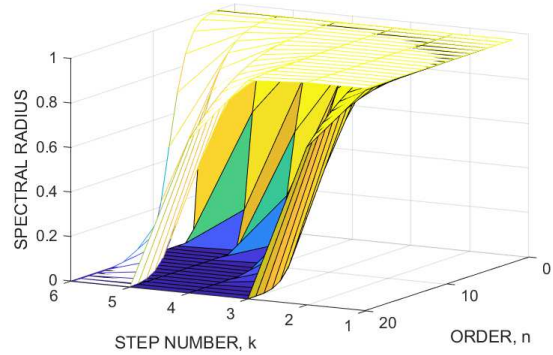,height=50mm}}
\begin{center}
\caption{\small{The Figure shows
the spectral radius  $ \ds \rho^{\ds h ~ n^2~ \frac{(k ~ n^{k+2} -(k-1)~ n^{k+1} - 2 ~ n^{k} - n + 2)}{(n-1)^2} } $ for Richardson iteration (for parameter estimation)
defined in (\ref{riem5}) (plotted as colored surface).
The spectral radius for Richardson iteration with improved convergence rate
described in \cite{st19},  $\ds \rho^{ \ds 2 h \{ \frac{n^{k+3} - n^4}{(n-1)^3}  - (k-1) \{ \frac{n^3}{(n-1)^2} + \frac{k}{2(n-1)} \}  + k (n+2) \} }$ is  plotted as white surface.
Both surfaces are plotted for  the spectral radius of ill-conditioned case
as functions of the order $n$ and step number $k$ for $h=1$.
 }}
\label{figsurface2}
\end{center}
\end{figure}

\begin{figure}
\centerline{\psfig{figure=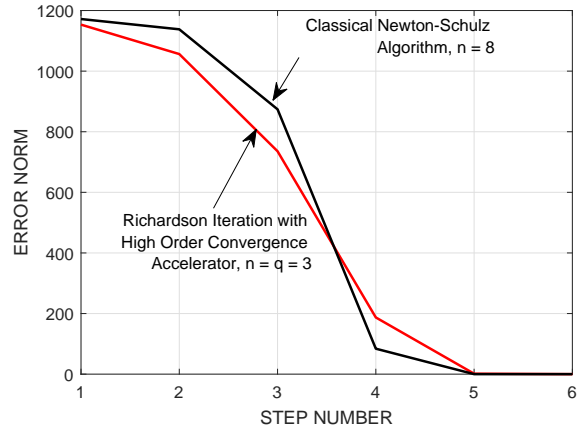,height=60mm}}
\begin{center}
\caption{\small{The Figure shows
the error of the estimated parameters of Richardson iteration with convergence accelerator of the third order (plotted with a red line) and the parameter error
of classical eight-order Newton-Schulz algorithm, plotted with a black line.
Both algorithms converge in five steps.}}
\label{figcomp1}
\end{center}
\end{figure}

\noindent
Numerical calculation of the parameter vector $\tt_*$
for the system (\ref{ae1}) where the ill-conditioned
SPD information matrix $A$  associated with
the system with harmonic regressor with three frequencies,
\cite{st15}, \cite{st10}, \cite{st14} is chosen for comparisons.
The performance evaluation is presented in the following three parts.
\\ 1) The convergence rate of new matrix inversion algorithm (\ref{gamk11}) - (\ref{ermoex1})
is compared to the convergence rate of recent algorithm
with improved convergence rate
described in \cite{st19}, see Figure~\ref{figsurface3}.
The Figure shows that convergence rate improvements
are more pronounced for higher orders and larger step numbers.
\\ 2)
Comparison of the convergence rate of the parameter estimation
algorithm (\ref{ttk1}) and the Richardson iteration with improved convergence rate
described in \cite{st19},  is presented in Figure~\ref{figsurface2}.
The algorithm (\ref{ttk1}) with high order convergence accelerator
improves essentially the convergence rate compared to existing algorithms
even for lower orders and small step numbers.
Indeed, comparison of the  Figure~\ref{figsurface3} and
Figure~\ref{figsurface2} shows that new algorithms are the most beneficial
in  the Richardson framework.
\\ 3) Finally, the performance evaluation of the algorithm (\ref{ttk1}) with respect to
classical Newton-Schulz algorithm is presented in  Figure~\ref{figcomp1}.
The Figure shows that the convergence rate of the
Richardson iteration with convergence accelerator of the order three is
comparable to the rate of classical Newton-Schulz algorithm (applied to the parameter estimation problem) of the order eight. Parameter estimation accuracy of the Richardson iteration (\ref{ttk1})
is about five times higher compared to the accuracy of classical Newton-Schulz algorithm in finite digit
calculations.

\section{Conclusion}
\label{conc}
\noindent
This paper shows that the most general and well-known Newton-Schulz iteration is fast power series expansion and presents unified framework and tool-kit for power series factorization and reduction of the computational complexity.
The framework allows reduction of complexity of many algorithms and factorization of the algorithm of the order $45$ that requires $10$ mmm only is presented as example.
\\
Main result of the paper is   new composite power series expansion for
Newton-Schulz iteration with high degree of parallelism for the
convergence rate improvement and computational efficiency. Comparative
analysis of the convergence rates of new algorithms and
exiting ones is performed via explicit transient models.
New algorithms have faster convergence  than known
Newton-Schulz iterations.
Moreover,  new expansion resulted in significant
improvement of the convergence rate of Richardson iteration
for which recursive and computationally efficient version was developed.
The results were also confirmed by simulations.
\\ The paper opens new opportunities for convergence rate improvement of
Newton-Schulz and Richardson iterations
via computationally efficient composite expansions
to be implemented on parallel machines with different computational performance.



\newpage

\section{Appendix. Factorization Tool-Kit: A Unified Approach}
\label{appen}
Some of known factorizations of the Newton-Schulz iteration are presented
in the following unified framework (see the Tables~\ref{t1} - \ref{t3} below):

\baq
  Z &=& \{ \sum_{j=0}^{h-1} Y^j \} ~ X
   \label{z1a} \\
  Y &=& I - X A \label{xanew} \\
  Z &=& \{ \sum_{j=0}^{w-1} Y^{(p+1) j} \} ~ \{ \sum_{d=0}^{p} Y^d \} ~ X    \label{ztab1} \\
   h &=& w ~ ( p + 1), ~~ p = 0,1,2, ..., ~ w = 1,2,3, ...  \label{hape} \\
  Z_1 &=& (I + Y \{ \sum_{j=0}^{w-1} Y^{(p+1) j} \} ~ \{ \sum_{d=0}^{p} Y^d \}) ~ X  \label{ztab2}  \\
   h_1 &=& h + 1 \label{hape1}
\eaq
where equations (\ref{z1a}) and (\ref{xanew}) represent the Newton-Schulz iteration
and equations (\ref{ztab1}), (\ref{ztab2})  represent the factorizations
$Z$ and $Z_1$  of orders  $h$ is $h_1$ respectively.
\\
The factorization (\ref{ztab1}) is valid for the orders, which are presented as
the composite numbers\footnote{A composite number is a positive integer that has at least one divisor other than one and itself or can be formed by multiplying two smaller positive integers},~ $ h = 2,4,6,8,9,10,12,14,15,16,18$
(excepting $h=2$).
For the orders which represent the prime numbers\footnote{A prime number is a positive integer that has exactly two distinct whole number factors (or divisors), namely one and the number itself},~$ h_1 = 3,5,7,11,13,17,19$  the factorization (\ref{ztab2}) is valid. Notice that factorization for the order $h=10$ is presented in both forms following \cite{solor}, see Table~\ref{t2}.
\\ Notice that the factorization defined in (\ref{ztab1}) is not unique.
For example the Newton-Schulz iteration of order $h=18$ can be factorized in the following
four ways: a) $p=5$, $w = 3$, b) $p=2$, $w = 6$, c) $p=8$, $w = 2$, d) $p=1$, $w = 9$
which have different number of mmm in implementation.
\\ Notice also that the factorization  (\ref{ztab2}) is simple application of the idea
known as Schr\"oder–Traub sequence, \cite{sen1}, \cite{tra1} to the polynomial
factorized in (\ref{ztab1}). The factorization  (\ref{ztab2}) increases the order of (\ref{ztab1})
by one, see (\ref{hape}) and (\ref{hape1}).
\\ Moreover, nested application of the factorizations (\ref{ztab1}) and  (\ref{ztab2})
implies additional order reduction and improvement of the computational efficiency
(mainly for high orders), see for example nested factorization for
$h=15$ with $7$ mmm in Table~\ref{t3}.
\\ Finally, the Tables~\ref{t1} - \ref{t3} which end with the algorithm of the nineteenth order
can be easily extended for higher orders (for any order) using nested applications of the
factorizations (\ref{ztab1}), (\ref{ztab2}) and the tables of prime factors, \cite{wiki},
providing computationally efficient
and implementable solutions  for higher orders.
\\ However, the factorizations which require computational efforts and additional memory
may result in error accumulation in finite-digit calculations.
The factorizations in high order Newton-Schulz  iterations can be seen as the additional
sub-steps (nested calculations for order reduction).
The idea of order reduction is also associated with the Newton-Schulz iteration,
see Section~\ref{psf}.
Therefore Newton-Schulz  algorithms of low orders ($h = 2,3$) being iterated for a number of steps
can be applied instead of factorized Newton-Schulz  iterations of higher orders for the sake of
robustness and efficiency.
\\ Notice that two and three steps of the
second order Newton-Schulz iteration are equivalent to one step
of the fourth and eighth order iterations with $4$ and $6$ mmm respectively, see Table~\ref{t1}, $h=4$ and Table~\ref{t2}, $h = 8$.
Robustness, efficiency and accuracy arguments motivate application of the second order \cite{sch}- \cite{beniz} and the third order \cite{isa}  Newton-Schulz iterations  instead of higher orders
in some cases. However, application of the eleventh order algorithm, see Table~\ref{t2} and
\cite{stan11} requires also $6$ mmm and provides faster convergence
than three steps of the second order Newton-Schulz iteration.
\\ Therefore the proper choice of the order and factorization that is made for each particular
application should represent the trade-off between the robustness and convergence rate.
\newpage
\begin{table*}[htbp]
 \centering
 \begin{tabular}{|c | c | c | c |}
 \hline
             &             &        &    \\ [1.5ex]
 Order           &     Factorization         & Unified       & Ref.  \\ [1.5ex]
  $h$ &      & Factorization &   \\ [1.5ex]
 \hline\hline
      &     &  &   \\
 $h=2$ & $ \ds (I + Y)~X$  & $\ds \sum_{j=0}^{w-1} Y^{(p+1) j}
  \sum_{d=0}^{p} Y^d~ X$   & \cite{sch}- \cite{beniz} \\
    &     & $p=1$, $w = 1$ &   \\
 \hline
      &     &  &   \\
 $h=3$ & $ \ds (I + Y ( I + Y))~X$  & $\ds (I + Y \sum_{j=0}^{w-1} Y^{(p+1) j}
  \sum_{d=0}^{p} Y^d)~ X$   & \cite{sen1}, \cite{li1},\cite{pan}   \\
    &     & $p=1$, $w = 1$ &   \\
 \hline

      &     &  &   \\
 $h=4$ & $ \ds  (I + Y^2) (I + Y)~X$  & $\ds  \sum_{j=0}^{w-1} Y^{(p+1) j}
  \sum_{d=0}^{p} Y^d ~ X$   & \cite{esma}  \\
    &     & $p=1$, $w = 2$ &   \\
  \hline

       &     &  &   \\
 $h=5$ & $ \ds  (I + Y (I + Y^2) (I + Y))~X$  & $\ds (I + Y \sum_{j=0}^{w-1} Y^{(p+1) j}
  \sum_{d=0}^{p} Y^d )~ X$   & \cite{solor}  \\
    &     & $p=1$, $w = 2$ &   \\
       & $\ds = (I + Y + Y^2 + Y^2 (Y + Y^2))~X $    &  &   \\
         &     &  &   \\
       \hline
 $h=6$ & $ \ds  (I + Y^3) (I + Y  + Y^2 )~X$  & $\ds \sum_{j=0}^{w-1} Y^{(p+1) j}
  \sum_{d=0}^{p} Y^d ~ X$   &    \\
    &     & $p=2$, $w = 2$ &   \\
  \hline

     &     &  &   \\
 $h=7$ & $ \ds (I + (Y + Y^4) (I + Y + Y^2) )  ~X$  & $\ds  (I + Y  \sum_{j=0}^{w-1} Y^{(p+1) j}
  \sum_{d=0}^{p} Y^d) ~ X$   & \cite{sol7} \\
    &     & $p=2$, $w = 2$ &   \\ [1ex]
  \hline
\end{tabular}
\\
\caption{Factorization of the algorithms of orders $h=2,...,7$ } \label{t1}
\end{table*}


\begin{table*}[htbp]
 \centering
 \begin{tabular}{|c | c | c | c |}
 \hline
             &             &        &    \\ [1.5ex]
 Order           &      Factorization         & Unified       & Ref.  \\ [1.5ex]
  $h$ &      & Factorization &   \\ [1.5ex]
 \hline\hline
     &     &  &   \\
 $h=8$ & $ \ds  (I + Y^4) (I + Y^2) (I + Y)  X$  & $\ds  \sum_{j=0}^{w-1} Y^{(p+1) j}
  \sum_{d=0}^{p} Y^d X$   &  \\
    &     & $p=3$, $w = 2$ &   \\
 \hline

    &     &  &   \\
 $h=9$ & $ \ds  (I + Y^3 + Y^6) (I + Y  + Y^2 )~X$  & $\ds  \sum_{j=0}^{w-1} Y^{(p+1) j}
  \sum_{d=0}^{p} Y^d   X$   &  \\
    &     &   &   \\
  &  $ = (I + (I + Y^4) (I + Y^2) (Y + Y^2) ) ~ X $ &  $p=2$, $w = 3$  & \\
   &     &  &   \\
 \hline
  &     &  &   \\
 $h=10$ & $ \ds (I + Y (I + Y^3 + Y^6)  $   & $\ds (I + Y  \sum_{j=0}^{w-1} Y^{(p+1) j}
  \sum_{d=0}^{p} Y^d) X$   & \cite{solor}, \cite{jeb} \\
    & $\ds   (I + Y + Y^2)  )  ~ X = $     &     &   \\
    &     & $p=2$, $w = 3$ &   \\
     &  $  (I + Y^5) (I + (Y  + Y^2)(I + Y^2) ) ~ X $  &  $\ds = \sum_{j=0}^{w-1} Y^{(p+1) j}
  \sum_{d=0}^{p} Y^d X$  &   \\
   &     & $p=4$, $w = 2$   &   \\
\hline
&     &  &   \\
 $h=11$ & $ \ds  (I + Y (I + Y^5) (I + (Y  + Y^2) $   & $\ds  (I + Y  \sum_{j=0}^{w-1} Y^{(p+1) j}
  \sum_{d=0}^{p} Y^d) X$   & \cite{solor} \\
    & $\ds  (I + Y^2)  )  ~ X$     &     &   \\
    &     & $p=4$, $w = 2$ &   \\
   &  $ \ds  (I + Y (I  + (Y^2 + Y^4) (I + Y^4) ) $    & $  \ds  (I + Y  \sum_{j=0}^{w-1} Y^{(p+1) j}
  \sum_{d=0}^{p} Y^d) X $   &  \cite{stan11}  \\
    &  $\ds (I + Y)  ) ~ X  $  & $p=1$, $w = 5$ &   \\
 \hline
&     &  &   \\
 $h=12$ & $ \ds (I + Y^4 + Y^8)  $   & $\ds  \sum_{j=0}^{w-1} Y^{(p+1) j}
  \sum_{d=0}^{p} Y^d ~ X$   & \cite{solor} \\
    & $\ds  (I + Y + Y^2 + Y^3 )  ~ X$     &     &   \\
    &     & $p=3$, $w = 3$ &   \\
\hline
 &     &  &   \\
 $h=13$ & $ \ds (I + Y  (I + Y^4 + Y^8)  $   & $\ds  (I +  Y  \sum_{j=0}^{w-1} Y^{(p+1) j}
  \sum_{d=0}^{p} Y^d)~ X$   & \cite{solor} \\
    & $\ds  (I + Y + Y^2 + Y^3 )   )  ~ X$     &     &   \\
    &     & $p=3$, $w = 3$ &   \\ [1ex]
  \hline
\end{tabular}
\\
\caption{Factorization of the algorithms of orders $h=8,...,13$ } \label{t2}
\end{table*}
\newpage
\begin{table*}[htbp]
 \centering
 \begin{tabular}{|c | c | c | c |}
 \hline
             &             &        &    \\ [1.5ex]
 Order           &      Factorization         & Unified       & Ref.  \\ [1.5ex]
  $h$ &      & Factorization &   \\ [1.5ex]
 \hline\hline
 &     &  &   \\
  $h=14$   &  $ (I + Y^7) $   &  $\ds  \sum_{j=0}^{w-1} Y^{(p+1) j} \sum_{d=0}^{p} Y^d ~ X$    & \cite{solor}  \\
 &     &  &   \\
 & $ (I +  Y + Y^2 + Y^3 + Y^4 + Y^5 + Y^6) ~ X $    & $p=6$, $w = 2$ &   \\
 &     &  &   \\
  \hline
&     &  &   \\
 $h=15$ & $ \ds (I + Y^3 + Y^6 + Y^9 + Y^{12} ) $   & $\ds  \sum_{j=0}^{w-1} Y^{(p+1) j}
  \sum_{d=0}^{p} Y^d ~ X$   & \cite{solor} \\
    & $\ds  (I + Y + Y^2)    ~ X$     &     &   \\
    &     & $p=2$, $w = 5$ &   \\
    &   $ \ds = ( I + (I + (Y^3)^2)~( (Y^3)^2 + Y^3  ) ) $       &  &   \\
    &     &  &   \\
    &  $\ds  (I + Y + Y^2)    ~ X$       &  &   \\
    &     &  &   \\
   \hline
&     &  &   \\
 $h=16$ & $ \ds (I + Y^4 + Y^8 + Y^{12} )  $   & $\ds   \sum_{j=0}^{w-1} Y^{(p+1) j}
  \sum_{d=0}^{p} Y^d ~ X$   & \cite{solor} \\
    & $\ds   (I + Y + Y^2 + Y^3)  ~ X$     &     &   \\
    &     & $p=3$, $w = 4$ &   \\
\hline
 &     &  &   \\
 $h=17$ & $ \ds  (I + (Y + Y^2 + Y^3 + Y^4 ) $   & $\ds  (I +  Y  \sum_{j=0}^{w-1} Y^{(p+1) j}
  \sum_{d=0}^{p} Y^d)~ X$   & \cite{solor} \\
    & $\ds  (I + Y^4 + Y^8 + Y^{12} ) )  ~ X$     &     &   \\
    &     & $p=3$, $w = 4$ &   \\
  \hline
&     &  &   \\
 $h=18$ & $ \ds  (I + Y^6 + Y^{12} )  $   & $\ds   \sum_{j=0}^{w-1} Y^{(p+1) j}
  \sum_{d=0}^{p} Y^d ~ X$   & \cite{solor} \\
    & $\ds   (I + Y + Y^2 + Y^3 + Y^4 + Y^5)     ~ X$     &     &   \\
    &     & $p=5$, $w = 3$ &   \\
\hline

 &     &  &   \\
 $h=19$ & $ \ds  (I + (Y + Y^2)(I + Y^2 + Y^4 ) $   & $\ds  (I +  Y  \sum_{j=0}^{w-1} Y^{(p+1) j}
  \sum_{d=0}^{p} Y^d)~ X$   & \cite{solor} \\
    & $\ds  (I + Y^6 + Y^{12}) )   ~ X$     &     &   \\
    &     & $p=5$, $w = 3$ &   \\ [1ex]
 \hline
\end{tabular}
\\
\caption{Factorization of the algorithms of orders $h=14,...,19$ } \label{t3}
\end{table*}

\end{document}